\documentclass[12pt]{amsart}

\usepackage[toc]{appendix} 
\usepackage{latexsym}
\usepackage{amsfonts}
\usepackage{amssymb}
\usepackage{amsmath}
\usepackage{mathrsfs}
\usepackage{bbm}
\usepackage{dsfont}
\usepackage{hyperref}
\allowdisplaybreaks

\usepackage[toc]{appendix}

\newcommand{\be}{\begin{equation}}
\newcommand{\ee}{\end{equation}}
\newcommand{\bea}{\begin{eqnarray}}
\newcommand{\eea}{\end{eqnarray}}
\newcommand{\bean}{\begin{eqnarray*}}
\newcommand{\eean}{\end{eqnarray*}}
\newcommand{\brray}{\begin{array}}
\newcommand{\erray}{\end{array}}

\newtheorem{dfn}{Definition}[section]
\newtheorem{thm}[dfn]{Theorem}
\newtheorem{lmma}[dfn]{Lemma}
\newtheorem{ppsn}[dfn]{Proposition}
\newtheorem{crlre}[dfn]{Corollary}
\newtheorem{xmpl}[dfn]{Example}
\newtheorem{rmrk}[dfn]{Remark}

\newcommand{\bdfn}{\begin{dfn}\rm}
\newcommand{\bthm}{\begin{thm}}
\newcommand{\blmma}{\begin{lmma}}
\newcommand{\bppsn}{\begin{ppsn}}
\newcommand{\bcrlre}{\begin{crlre}}
\newcommand{\bxmpl}{\begin{xmpl}}
\newcommand{\brmrk}{\begin{rmrk}\rm}

\newcommand{\edfn}{\end{dfn}}
\newcommand{\ethm}{\end{thm}}
\newcommand{\elmma}{\end{lmma}}
\newcommand{\eppsn}{\end{ppsn}}
\newcommand{\ecrlre}{\end{crlre}}
\newcommand{\exmpl}{\end{xmpl}}
\newcommand{\ermrk}{\end{rmrk}}

\newcommand{\bbc}{\mathbb{C}}
\newcommand{\bbz}{\mathbb{Z}}

\newcommand{\bbn}{\mathbb{N}}
\newcommand{\bbr}{\mathbb{R}}

\makeatletter
\let\@wraptoccontribs\wraptoccontribs
\makeatother

\author{ S. Sundar}
\title{Representations of the weak Weyl commutation relation}

\begin{document}
\maketitle
\begin{abstract}
Let $G$ be a locally compact abelian group with Pontraygin dual $\widehat{G}$. Suppose $P$ is a closed subsemigroup of $G$ containing the identity element $0$. We assume that $P$ has dense interior and $P$ generates $G$. Let $U:=\{U_{\chi}\}_{\chi \in \widehat{G}}$ be a strongly continuous group of unitaries and let $V:=\{V_{a}\}_{a \in P}$ be a strongly continuous semigroup of isometries. We call  $(U,V)$ a weak Weyl pair if 
\[
U_{\chi}V_{a}=\chi(a)V_{a}U_{\chi}\]
for every $\chi \in \widehat{G}$ and for every $a \in P$. 

We work out the representation theory (the factorial and the irreducible representations) of the above commutation relation under the assumption that $\{V_{a}V_{a}^{*}:a \in P\}$ is a commuting family of projections. Not only does this generalise the results of \cite{Bracci1} and \cite{Bracci2}, our proof brings out the Morita equivalence that lies behind the results.  For $P=\bbr_{+}^{2}$, we  demonstrate that if we drop the commutativity assumption on the range projections, then the representation theory of the weak Weyl commutation relation becomes very complicated. 
\end{abstract}

\noindent {\bf AMS Classification No. :} {Primary 46L05 ; Secondary 81S05.}  

\section{Introduction}

The classical Stone-von Neumann theorem that asserts the uniqueness of the Weyl commutation relation
\[
U_sV_t=e^{its}V_tU_s\]
where $\{U_s\}_{s \geq 0}$ and $\{V_t\}_{t \geq 0}$ are strongly continuous $1$-parameter group of unitaries is a fundamental theorem in both quantum mechanics and in operator algebras. In \cite{Bracci1} and in \cite{Bracci2}, a weaker version of the above commutation relation is considered where it is assumed that $\{V_t\}_{t \geq 0}$ is only a semigroup of isometries. The representation theory (the factorial representations and the irreducible representations) of such relations was worked out by Bracci and Picasso in \cite{Bracci1} and in \cite{Bracci2}. Bracci and Picasso considers such a weak form of the commutation relation as the quantisation postulate for systems whose configuration space is semibounded like the half-line. This is because, on the half-line, though the position operator generates a group of unitaries, the  momentum operator (which is not self-adjoint)  generates only a semigroup of isometries.  

The purpose of this paper is twofold.   First,  we would like to bring out the $C^{*}$-algebraic reason for the validity of the results in \cite{Bracci1} and in \cite{Bracci2}. Secondly, we wish to extend slightly the resuts to systems with $d$ degrees of freedom
where $d \geq 2$. We work in the more general setting of subsemigroups of locally compact abelian groups. 

It is well known (\cite{Rosenberg}, \cite{Williams_Dana}) that the $C^{*}$-algebra that encodes the usual Weyl commutation relation  is Morita equivalent to $\bbc$. Stone-von Neumann theorem is then an immediate consequence of this Morita equivalence. 
We establish a similar reasoning here. We prove
that the $C^{*}$-algebra that encodes the weak version of the Weyl commutation relation considered in this paper is Morita equivalent to a commutative $C^{*}$-algebra. Thus, it follows at once that every factorial representation is a multiple of an irreducible representation and 
the irreducible representations are parameterised by the character space of the underlying commutative $C^{*}$-algebra. Moreover, in the one-dimensional case the commutative $C^{*}$-algebra is $C_{0}((-\infty,\infty])$. This provides a conceptual explanation for the 
results obtained in \cite{Bracci1} and in \cite{Bracci2}. 

The proof is based on the results obtained in \cite{Sundar_Ore} where a certain ``universal dynamical system" was constructed encoding all semigroups of isometries with commuting range projections. 
Strictly speaking, only half of the last statement was proved in \cite{Sundar_Ore} (see Remark \ref{Remark}) and we prove the other half in this paper modulo Morita equivalence.

The results obtained are next explained. 

Let $G$ be a locally compact, second countable,  Hausdorff abelian group. We denote the dual group of $G$ by $\widehat{G}$. We use additive notation for the group operations. Let $P \subset G$ be a closed subsemigroup containing $0$ such that $P-P=G$. Set $\Omega:=Int(P)$. We assume that $\Omega$ is dense in $P$. For $x,y \in G$, we write $x \leq y$ if $y-x \in P$ and $x<y$ if $y-x \in \Omega$

 Let $U:=\{U_\chi\}_{\chi \in \widehat{G}}$ be a strongly continuous group of unitaries and let $V:=\{V_{a}\}_{a \in P}$ be a strongly continuous semigroup of isometries. 
We call $(U,V)$ a weak Weyl pair if for every $\chi  \in \widehat{G}$ and for every $a \in P$, 
\[
U_\chi V_a=\chi(a)V_aU_\chi.\]

Let $(U,V)$ be a weak Weyl pair. For $a \in P$, let $E_a:=V_aV_{a}^{*}$. We say that $(U,V)$ has commuting range projections if  $\{E_{a}:a \in P\}$ is a commuting family of projections. Note that if $P=[0,\infty)$ or, more generally, if the preorder $\leq$ is a total order, every weak Weyl pair has commuting range projections. 

Examples of weak Weyl pairs with commuting range projections are given below. Let $A$ be a non-empty closed subset of $G$ which is $P$-invariant, i.e. $A+P \subset A$. Such subsets will be called $P$-spaces. Let $K$ be a Hilbert space whose dimension we denote by  $k$. Consider the Hilbert space $H:=L^{2}(A,K)$. 
For $\chi \in \widehat{G}$, let $U_\chi$ be the unitary on $H$ defined by 
\[
U_{\chi}f(y)=\chi(y)f(y).\]
Then, $U:=\{U_\chi\}_{\chi \in \widehat{G}}$ is a strongly continuous group of unitaries on $H$. 
For $a \in P$, let $V_{a}$ be the isometry on $H$ defined by 
\begin{equation*}
 \label{isometries}
V_{a}(f)(y):=\begin{cases}
 f(y-a)  & \mbox{ if
} y-a \in A,\cr
   &\cr
    0 &  \mbox{ if } y-a \notin A.
         \end{cases}
\end{equation*}
Then, $V=\{V_a\}_{a \in P}$ is a strongly continuous  semigroup of isometries on $H$. It is clear that $V$ has commuting range projections. It is routine to verify that $(U,V)$ is a weak Weyl pair. We call $(U,V)$ the weak Weyl pair associated to the $P$-space $A$ with multiplicity $k$. If we want to stress the dependence of $(U,V)$ on $A$ and $k$, we denote $(U,V)$ by $(U^{(A,k)},V^{(A,k)})$. 

The main theorem of this paper is stated below.
\begin{thm}
\label{main theorem}
We have the following. 
\begin{enumerate}
\item[(1)] Let $A$ be a $P$-space and let  $k \in \{1,2,\cdots\} \cup \{\infty\}$ be given. The weak Weyl pair $(U^{(A,k)},V^{(A,k)})$ is a factorial representation. Moreover, it is irreducible if and only if $k=1$. 
\item[(2)] Let $A,B$ be $P$-spaces and let $k,\ell \in \{1,2,\cdots \} \cup \{\infty\}$ be given. The weak Weyl pair $(U^{(A,k)},V^{(A,k)})$ is unitarily equivalent to $(U^{(B,\ell)},V^{(B,\ell)})$ if and only if $A=B$ and $k=\ell$. 
\item[(3)] Suppose $(U,V)$ is a weak Weyl pair with commuting range projections.  Assume that the von Neumann algebra generated by $\{U_{\chi},V_a: \chi \in \widehat{G},a \in P\}$ is a factor. Then, there exists a $P$-space $A$ and $k \in \{1,2,\cdots\} \cup \{\infty\}$ such that 
$(U,V)$ is unitarily equivalent to $(U^{(A,k)},V^{(A,k)})$. 

Thus, for weak Weyl pairs with commuting range projections, factorial representations are completely reducible. Moreover, irreducible weak Weyl pairs with commuting range projections are precisely those 
associated to $P$-spaces with multiplicity $1$. 
\end{enumerate}
\end{thm} 
For $P=[0,\infty)$, as already mentioned,  every weak Weyl pair has commuting range projections. Also, every $P$-space is either $\bbr$ or of the form $[a,\infty)$ for a unique $a \in \bbr$. It is now clear that the results of \cite{Bracci1} and \cite{Bracci2} for the semibounded case follow from Thm. \ref{main theorem}.

Morever, for irreducible weak Weyl pairs with commuting range projections, we have the following uniqueness result. We need a bit of notation. Let $U:=\{U_{\chi}\}_{\chi \in \widehat{G}}$ be a strongly continuous group of unitaries on a Hilbert space $H$. Then, $U$ determines a represenation $\pi_U$ of $ C_0(G) \cong C^{*}(\widehat{G})$ on $H$. We denote the unique closed subset of $G$ that corresponds to the ideal $Ker(\pi_U)$ by $Spec(U)$. 

\begin{crlre}
\label{uniqueness}
Let $(U,V)$ and $(\widetilde{U},\widetilde{V})$ be irreducible weak Weyl pairs with commuting range projections. Assume that $(U,V)$ acts on $H$ and $(\widetilde{U},\widetilde{V})$ acts on $\widetilde{H}$. Suppose $Spec(U)=Spec(\widetilde{U})$.
Then, there exists a unitary $X:H \to \widetilde{H}$ such that for $\chi \in \widehat{G}$ and $a \in P$, 
\[
XU_{\chi}X^{*}=\widetilde{U}_{\chi}~~;~~~XV_{a}X^{*}=\widetilde{V}_{a}.\]
\end{crlre}

What about  weak Weyl pairs which do not have commuting range projections ? For $P=\bbr_{+}^{2}$, we demonstrate that working out the irreducible weak Weyl pairs is a complicated task. We explain a procedure (preserving factoriality, type and irreducibility) that allows us to build weak Weyl pairs starting from a  representation of the free product $c_0(\bbn)\ast c_0(\bbn)$. We also prove that Corollary \ref{uniqueness}  no longer stays true if we relax the commutativity assumption on the range projections.

We end this introduction by mentioning that weak Weyl relations, for one degree of freedom, in the unbounded picture were analysed extensively in the literature. Some of the important papers that deal with the unbounded version are  \cite{Konrad}, \cite{Konrad1}, \cite{Jorgensen_Muhly}, \cite{Arai2}, \cite{Arai3}, and \cite{Arai1}. We do not touch the unbounded version here. The author is of the belief that the $C^{*}$-algebra machinery may not be sufficient to handle domain issues.

All the Hilbert spaces considered in this paper are assumed to be separable. Moreover, we use the convention that the inner product is linear in the first variable. 

\section{The equivalence between $Isom_c(P)$ and $Rep(C_0(Y_u)\rtimes G)$}

For the rest of this paper, $G$ stands for an arbitrary but a fixed locally compact, second countable, Hausdorff abelian group. The letter $P$ stands for a closed subsemigroup of $G$ containing the identity element $0$. We assume that $\Omega:=Int(P)$ is dense in $P$. We also assume $P-P=G$. We first review, from \cite{Sundar_Ore}, the construction of the universal dynamical system that encodes the isometric representations of $P$ with commuting range projections. 

Let $\mathcal{C}(G)$ be  the set of closed subsets of $G$ endowed with the Fell topology. Let 
\[
Y_u:=\{A \in \mathcal{C}(G): A \neq \emptyset, -P+A \subset A\}.\]
Endow $Y_u$ with the subspace topology inherited from the Fell topology on $\mathcal{C}(G)$. Then, $Y_u$ is a locally compact, second countable, Hausdorff space. Moreover, the map
\[
Y_u \times G \ni (A,x) \to A+x \in Y_u\]
defines an action of $G$ on $Y_u$. Set 
\[
X_u:=\{A \in Y_u: -P \subset A\}=\{A \in Y_u:0 \in A\}.\]
Then, $X_u$ is a compact subset of $Y_u$. Clearly, $X_u+P \subset X_u$. If $(s_n)$ is a cofinal sequence in $\Omega$, then 
\[Y_u=\bigcup_{n \geq 1}(X_u-s_n).\]
\textbf{Notation:} For $f \in C_{c}(G)$, let $\widetilde{f}:Y_u \to \bbc$ be defined by
\begin{equation}
\label{generating set}
\widetilde{f}(A):=\int f(x)1_{A}(x)dx.\end{equation}
Then, $\widetilde{f} \in C_{c}(Y_u)$. Moreover, $\{\widetilde{f}:f \in C_{c}(G)\}$ generates $C_0(Y_u)$. For $f \in L^1(G)$, we define $\widetilde{f} \in C_0(Y_u)$ exactly as in Eq. \ref{generating set}.

For the proof of the above assertions, we refer the reader to \cite{Hilgert_Neeb} and \cite{Sundar_Ore}. The reader is also recommended to consult Section 5 of \cite{Piyasa_Sundar}.

\begin{lmma}
\label{embedding}
The map $Y_u \ni A \to 1_{A} \in L^{\infty}(G)$ is a topological embedding. Here, $L^{\infty}(G)$ is identified with $L^{1}(G)^{*}$ and $L^{\infty}(G)$ is endowed with  the weak $^*$-topology. 
\end{lmma}
\textit{Proof.} The fact that $Y_u \ni A \to 1_{A} \in L^{\infty}(G)$ is a continuous injection follows from Prop. II.13\footnote{In \cite{Hilgert_Neeb}, it is assumed that $P$ is a Lie semigroup. However, the proof of Prop. II.13 given in \cite{Hilgert_Neeb} works for subsemigroups (with dense interior) of locally compact abelian groups.} of \cite{Hilgert_Neeb}. Suppose $(A_n)$ is a sequence in $Y_u$ and let $A \in Y_u$. Suppose $1_{A_n} \to 1_{A}$ in $L^{\infty}(G)$. Then, 
$\widetilde{f}(A_n) \to \widetilde{f}(A)$ for every $f \in C_{c}(G)$. Since $\{\widetilde{f}:f \in C_{c}(G)\}$ generates $C_0(Y_u)$, it follows that $A_n \to A$ in $Y_u$. Hence the proof. \hfill $\Box$

\begin{rmrk}
\label{Remark}
For a strongly continuous semigroup of isometries $V:=\{V_{a}\}_{a \in P}$, we say $V$ has commuting range projections if $\{V_{a}V_{a}^{*}:a \in P\}$ is a commuting family of projections. We also call a strongly continuous semigroup of isometries indexed by $P$
 an isometric representation of $P$. 

Let $Isom_c(P)$ be the collection (up to unitary equivalence) of isometric representations of $P$ with commuting range projections. Let $\mathcal{G}_u:=X_u \rtimes P$ be the Deaconu-Renault groupoid considered in \cite{Sundar_Ore}. 
Let $Rep(C^*(\mathcal{G}_u))$ be the collection (up to unitary equivalence) of non-degenerate representations of $C^{*}(\mathcal{G}_u)$. It follows from Theorem 7.4 of \cite{Sundar_Ore} that there exists an injective map 
\[
Isom_c(P) \ni V \to \pi_V \in Rep(C^{*}(\mathcal{G}_u)).\]
However, it was not proved in \cite{Sundar_Ore} that the above map is surjective.  More precisely, the inverse map was not constructed in \cite{Sundar_Ore}. In this paper, we correct this deficiency modulo Morita equivalence by passing to the transformation groupoid $Y_u \rtimes G$ which is equivalent to $\mathcal{G}_u$.

We show in this paper how to construct maps $\Phi:Rep(C_0(Y_u)\rtimes G) \to Isom_c(P)$ and $\Psi:Isom_c(P) \to Rep(C_0(Y_u)\rtimes G)$ which are inverses of each other.  To avoid repetition, we will be economical with details and omit proofs that require minor modifications of the arguments presented in \cite{Sundar_Ore}. 

The reason for preferring the transformation groupoid  $Y_u\rtimes G$  over the Deaconu-Renault groupoid  $X_u \rtimes P$ is that   Takai duality is readily available for the transformation groupoid $Y_u \rtimes G$. 
\end{rmrk}

Let $(\pi,W)$ be a covariant representation of the dynamical system $(C_0(Y_u),G)$ on a separable Hilbert space $K$.   Denote the algebra of bounded Borel functions on $Y_u$ by $B(Y_u)$. For $\phi \in B(Y_u)$ and $x \in G$, let 
$L_x\phi \in B(Y_u)$ be defined by $L_x(\phi)(A)=\phi(A-x)$. Denote the extension of $\pi$ to $B(Y_u)$, obtained via the measurable functional calculus, by $\pi$ itself. Then, we have the following covariance relation: for $\phi \in B(Y_u)$ and $x \in G$, 
\[
W_x\pi(\phi)W_x^{*}=\pi(L_x(\phi)).\] 

Set $H:=\pi(1_{X_u})K$. Since $X_u+P \subset X_u$, it follows that the subspace $H$ is invariant under $\{W_a:a \in P\}$. For $a \in P$, let $V_{a}$ be the operator on $H$ defined by $V_{a}:=W_{a}|_{H}$. Then, $V:=\{V_{a}\}_{a \in P}$ is a strongly continuous semigroup of isometries. Moreover, the collection $\{V_{a}V_{a}^{*}:a \in P\}$ is a commuting family of projections on $H$. If we want to stress the dependence of $V$ on $(\pi,W)$, we denote $V$ by $V^{(\pi,W)}$.

Thus, we get a map denoted $\Phi$ 
\[
Rep(C_0(Y_u) \rtimes G) \ni (\pi,W) \to V^{(\pi,W)} \in Isom_c(P).\]
We next explain how to construct the inverse of  the map $\Phi$. 

Let $V:=\{V_{a}\}_{a \in P}$ be a strongly continuous semigroup of isometries acting on a Hilbert space $H$ with commuting range projections. Let $(W,K)$ be the minimal unitary dilation of $V$. This means the following.
\begin{enumerate}
\item[(1)] The Hilbert space $K$ contains $H$ as a closed subspace.
\item[(2)] $W=\{W_x\}_{x \in G}$ is a strongly continuous group of unitaries on $K$.
\item[(3)] For $a \in P$ and $\xi \in H$, $W_a\xi=V_a\xi$.
\item[(4)] The union $\bigcup_{a \in P}W_{a}^{*}H$ is dense in $K$.
\end{enumerate}

For $x \in G$, let $E_x$ be the projection onto the subspace $W_xH$. Observe the following. 
\begin{enumerate}
\item[(1)] For $x, y\in G$, $W_xE_yW_{x}^{*}=E_{x+y}$. 
\item[(2)] Since $W_aH \subset H$ for $a \in P$, it follows that for $x,y \in G$ with $x \leq y$, $E_{x} \geq E_y$. 
\item[(3)] Let $x, y \in G$ be given. Then, $E_xE_y=E_yE_x$. If $x,y \in P$, this follows from the fact that $V$ has commuting range projections. Note that $\displaystyle \bigcup_{c \in \Omega}(\Omega-c)=\Omega-\Omega=G$.  Making use of the fact that $(\Omega-c_1)\cup (\Omega-c_2) \subset \Omega-(c_1+c_2)$ for $c_1,c_2 \in \Omega$,  choose $a,b,c \in \Omega$ such that $x=a-c$ and $y=b-c$. Then, by $(1)$, $E_{x}=W_{c}^{*}E_{a}W_{c}$ and $E_{y}=W_{c}^{*}E_{b}W_{c}$.  Since $E_a$ and $E_b$ commute, $E_x$ and $E_y$ commute. 
\item[(4)] The map $G \ni x \to E_{x}=W_{x}E_0W_{x}^{*} \in B(K)$ is strongly continuous. 
\end{enumerate}

Let $\mathcal{D}$ be the $C^{*}$-subalgebra of $B(K)$ generated by $\Big\{ \int f(x)E_xdx: f \in L^1(G)\Big\}$. Let $\chi$ be a character of $\mathcal{D}$. Arguing exactly as in  the proof of Prop. 4.3 of \cite{Sundar_Ore} and as in the proof of Prop. 4.6 of \cite{Sundar_Ore}, we see that there exists a unique element in $Y_u$, denoted $A_{\chi}$, such that 
\[
\chi\Big(\int f(x)E_xdx\Big)=\int f(x)1_{A_\chi}(x)dx\]
for every $f \in L^{1}(G)$. 
It follows from the above equality and Lemma \ref{embedding} that the map $\widehat{\mathcal{D}} \ni \chi \to A_{\chi} \in Y_u$ is a topological embedding. Via this embedding, we view $\widehat{\mathcal{D}}$ as a subspace of $Y_u$. 

\begin{lmma}
The subset $\widehat{\mathcal{D}}$ is a closed subset of $Y_u$. 
\end{lmma}
\textit{Proof.} Let $(\chi_n)$ be a sequence in $\widehat{\mathcal{D}}$ such that $A_{\chi_n} \to A$ for some $A \in Y_u$. Passing to a subsequence, if necessary, we can assume that $(\chi_n) \to \chi$ for some element in $\mathcal{D}^{*}$ where $\mathcal{D}^{*}$, the dual of $\mathcal{D}$,  is given the weak $^*$-topology. We claim that $\chi \neq 0$. 

Choose $f \in C_{c}(G)$ such that $\int f(x)1_{A}(x)dx=1$. Then, 
\begin{align*}
\chi\big(\int f(x)E_xdx\big)&=\lim_{n \to \infty}\chi_{n}\big(\int f(x)E_xdx \big)\\&=\lim_{n \to \infty}\int f(x)1_{A_{\chi_n}}(x)dx\\& = \int f(x)1_{A}(x)dx\\&=1.\end{align*}
Thus, $\chi$ is non-zero. This proves the claim. 

Since $\chi \neq 0$, $\chi$ is a character of $\mathcal{D}$. The fact that $\chi$ is a character and the continuity of the map $\widehat{\mathcal{D}} \ni \psi \to A_{\psi} \in Y_u$ implies that $A_{\chi}=\lim_{n \to \infty}A_{\chi_n}=A$. Hence, $\widehat{\mathcal{D}}$ is a closed subset of $Y_u$. Hence the proof. \hfill $\Box$

Let $Res:C_{0}(Y_u) \to C_{0}(\widehat{\mathcal{D}})$ be the restriction map and let $G:C_{0}(\widehat{\mathcal{D}}) \to \mathcal{D}$ be the inverse of the Gelfand transform. Define $\pi:C_{0}(Y_u) \to \mathcal{D} \subset B(K)$ by 
\[
\pi(\phi)=G\circ Res(\phi).\]

Note that $\pi$ is the unique $^*$-homomorphism such that 
\begin{equation}
\label{defining equality}
\pi(\widetilde{f})=\int f(x)E_xdx
\end{equation}
for $f \in L^1(G)$. 

Using Eq. \ref{defining equality}, the equality $W_xE_yW_{x}^{*}=E_{x+y}$ and the fact that $\{\widetilde{f}: f \in C_{c}(G)\}$ generates $C_0(Y_u)$, it is routine to verify that $(\pi,W)$ is a covariant representation of the dynamical system 
$(C_0(Y_u),G)$. 
We denote the extension of $\pi$ to $B(Y_u)$ obtained via the Borel functional calculus by $\pi$ itself. Then, 
\[
W_x\pi(\phi)W_{x}^{*}=\pi(L_x(\phi))\]
for $x \in G$ and $\phi \in B(Y_u)$. 

We record below a few elementary properties of the representation $\pi$. 

\begin{lmma}
\label{projection}
For $x \in G$, $\pi(1_{X_u+x})=E_x$. In particular, $\pi(1_{X_u})$ is the orthogonal projection onto $H$.
\end{lmma}
\textit{Proof.} Thanks to the covariance relation, it suffices to prove that $\pi(1_{X_u})=E_0$. Let $(O_n)$ be a decreasing sequence of open sets with compact closure such that $\{O_n: n\geq 1\}$ forms a neighbourhood base
at $0$. Set $E_n:=O_n \cap -P$ and let $f_n:=\frac{1}{\mu(E_n)}1_{E_n}$. Here, $\mu$ is the Haar measure on $G$. 
We leave it to the reader to verify that $\widetilde{f_n} \to 1_{X_u}$ pointwise. Moreover, $\widetilde{f_n}$ is uniformly bounded. 

Thus, 
\begin{align*}
\pi(1_{X_u})&=\lim_{n \to \infty}\pi(\widetilde{f_n})\\
&=\lim_{n \to \infty}\int f_n(x)E_xdx\\
&=E_0.
\end{align*}
In the above chain of equalities, the limit is to be understood in the strong operator topology. Hence the proof. \hfill $\Box$

\begin{ppsn}
The representation $\pi$ is non-degenerate. 
\end{ppsn}
\textit{Proof.} Let $\xi \in K$ be such that $\langle \pi(\phi)\eta|\xi\rangle=0$ for every $\phi \in C_0(Y_u)$ and for every $\eta \in K$. Taking $\phi=\widetilde{f}$ for $f \in C_c(G)$, we see that 
\[
\langle \pi(\widetilde{f})\eta|\xi\rangle=\int f(x)\langle E_x\eta|\xi\rangle dx=0.\]
As the above equality happens for every $f \in C_{c}(G)$ and the map $G \ni x \to E_x \in B(K)$ is strongly continuous, we deduce that for $x \in G$ and $\eta \in K$, 
$\langle E_x\eta|\xi\rangle=0$. Thus, $\xi$ is orthogonal to $W_{a}^{*}H$ for every $a \in P$. But the union $\bigcup_{a \in P}W_{a}^{*}H$ is dense in $K$. Therefore, $\xi=0$. This completes the proof. \hfill $\Box$

To denote the dependence of $(\pi,W)$ on $V$, we denote $(\pi,W)$ by $(\pi^V,W^V)$. This way, we obtain a map, denoted $\Psi$, 
\[
Isom_{c}(P) \ni V \to (\pi^V,W^V) \in Rep(C_0(Y_u)\rtimes G).\]

\begin{thm}
\label{equivalence}
The map \[\Phi:Rep(C_0(Y_u) \rtimes G) \ni (\pi,W) \to V^{(\pi,W)} \in Isom_{c}(P)\] and the map \[\Psi:Isom_{c}(P) \ni V \to (\pi^V,W^V) \in Rep(C_0(Y_u)\rtimes G)\] are inverses of each other. 
\end{thm}
\textit{Proof.} Let $(\pi,W) \in Rep(C_0(Y_u)\rtimes G)$ be given. Suppose that $(\pi,W)$ acts on $K$. Set $V:=V^{(\pi,W)}$. We need to show that $(\pi^V,W^V)=(\pi,W)$. 

First, we claim that $W$ is the minimal unitary dilation of $V$. Recall that $H=\pi(1_{X_u})K$ and $V=\{V_{a}\}_{a \in P}$ is the restriction of $\{W_{a}\}_{a \in P}$ onto $H$. 
Thus, $W$ is a dilation of $V$. It is enough to show that $\bigcup_{a \in P}W_{a}^{*}H$ is dense in $K$. 

Let $(s_n)$ be an increasing cofinal sequence in $\Omega$.  Observe that $1_{X_u-s_n} \nearrow 1_{Y_u}$. Thus, $\pi(1_{X_u-s_n})=W_{s_n}^{*}\pi(1_{X_u})W_{s_n} \nearrow 1$ strongly. Clearly, $\pi(1_{X_u-{s_n}})$ is the orthogonal projection onto $W_{s_n}^{*}H$. 
Thus, $\bigcup_{n \geq 1}W_{s_n}^{*}H$ is dense in $K$. This proves that $W$ is the minimal unitary dilation of $V$. Thus, $W^{V}=W$. 

Next, we show that $\pi=\pi^V$. By definition, $\pi(1_{X_u+x})$ is the orthogonal projection onto $W_xH$ and by Lemma \ref{projection}, $\pi^V(1_{X_u+x})=E_x$ which is the orthogonal projection onto $W_xH$. Thus, 
\begin{equation}
\label{basic}
\pi(1_{X_u+x})=\pi^V(1_{X_u+x})
\end{equation}
for every $x \in G$. 

For a compact subset $F$ of $G$ and for an open subset $O$ of $G$, define
\begin{align*}
\mathcal{U}_{F}:&=\{A \in Y_u: A \cap F=\emptyset\},\\
\mathcal{U}_{O}:&=\{A \in Y_u: A \cap O\neq \emptyset\},\\
\mathcal{U}^{'}_{O}:&=\{A \in Y_u: A \cap O=\emptyset\}.
\end{align*}
The sets $U_{F} \cap \mathcal{U}_{O_1} \cap \mathcal{U}_{O_2} \cap \cdots \mathcal{U}_{O_n}$, as $F$ and $O_i$'s vary, form a basis for the Fell topology on $Y_u$. Thus, it suffices to show that for every compact set $F$ and for every open 
set $O$, $\pi(1_{\mathcal{U}_F})=\pi^V(1_{\mathcal{U}_F})$ and $\pi(1_{\mathcal{U}_{O}})=\pi^V(1_{\mathcal{U}_O})$. 

Fix an open set $O$ of $G$. Let $D:=\{x_1,x_2,\cdots\}$ be a dense subset of $O$. Let $A \in Y_u$ be given. Observe $A+\Omega=\bigcup_{a \in A}(a+\Omega)$ is an open set contained in $A$. Thus, $A+\Omega \subset Int(A)$. Since $0 \in \overline{\Omega}$, $Int(A)$ is dense in $A$.  Thus, for $A \in Y_u$, $A \cap O \neq \emptyset$ if and only if $Int(A) \cap O\neq \emptyset$ if and only if $A \cap D\neq \emptyset$. 
Therefore, for $A \in Y_u$, 
\[
1_{\mathcal{U}_O}(A)=\sup_{n \geq 1}1_{A}(x_n)=\sup_{n \geq 1}1_{X_u+x_n}(A).\]

By Eq. \ref{basic} and by Borel functional calculus, we obtain $\pi(1_{\mathcal{U}_O})=\pi^V(1_{\mathcal{U}_O})$. Since $\mathcal{U}^{'}_{O}$ is the complement of $\mathcal{U}_{O}$, it follows that 
\begin{equation}
\label{basic1}
\pi(1_{\mathcal{U}^{'}_{O}})=\pi^{V}(1_{\mathcal{U}^{'}_{O}})
\end{equation}
for every open set $O$ of $G$. 

Let $F$ be a compact subset of $G$. Choose a decreasing sequence of open sets $(O_n)$ such that $\{O_n: n \geq 1\}$ forms a base at $F$. This means that if $O$ is an open set that contains $F$, then
$O_n \subset O$ eventually. Note that for a closed subset $A$ of $G$, $A \cap F=\emptyset$ if and only if $A \cap O_n=\emptyset$ eventually. Thus, 
\[
1_{\mathcal{U}_F}=\limsup_{n \to \infty}1_{\mathcal{U}^{'}_{O_n}}.\]
By Borel functional calculus and by Eq. \ref{basic1}, we have $\pi(1_{\mathcal{U}_F})=\pi^V(1_{\mathcal{U}_F})$. Hence, $\pi=\pi^V$. This completes the proof of the assertion $\Psi\circ \Phi=Id$. 

Let $V \in Isom_{c}(P)$ be given. Suppose that $V$ acts on $H$.   Set $(\pi,W)=(\pi^V,W^V)$ and let $K$ be the Hilbert space on which $(\pi,W)$ acts. Let $\widetilde{V}=V^{(\pi,W)}$. For $a \in P$, $V_{a}$ is the restriction of $W_{a}$ to $H$ and $\widetilde{V}_a$ is the restriction of $W_{a}$ to $\pi(1_{X_u})H$. By Lemma \ref{projection},
we have $\pi(1_{X_u})H=H$. Consequently, $\widetilde{V}=V$. Hence $\Phi \circ \Psi=Id$. The proof is now complete. \hfill $\Box$

\section{Proof of the main theorem}
With Thm. \ref{equivalence} in hand, the conceptual explanation for Thm. \ref{main theorem} is quite simple. Having a unitary group, indexed by $\widehat{G}$, implementing the Weyl commutation relation
is equivalent to having a unitary group implementing the dual action on $C_0(Y_u)\rtimes G$. Then, Thm. \ref{main theorem} is a straightforward consequence of Takai duality. We explain some details below. 

Let $\mathcal{W}_c(P,\widehat{G})$ denote the collection (up to unitary equivalence) of weak Weyl pairs with commuting range projections. Consider the dual action of $\widehat{G}$ on $C_0(Y_u) \rtimes G$.
We prove below that $\mathcal{W}_{c}(P,\widehat{G}) \cong Rep((C_0(Y_u)\rtimes G)\rtimes \widehat{G})$. 

\begin{thm}
\label{equivalence1}
There exist maps \[\Psi:\mathcal{W}_{c}(P,\widehat{G}) \to Rep((C_0(Y_u)\rtimes G)\rtimes \widehat{G})\] and \[\Phi:Rep((C_0(Y_u)\rtimes G)\rtimes \widehat{G})\to \mathcal{W}_{c}(P,\widehat{G})\]
such that $\Phi$ and $\Psi$ are inverses of each other. 
\end{thm}
\textit{Proof.} Let $((\pi,W),U) \in Rep((C_0(Y_u)\rtimes G) \rtimes \widehat{G})$. Suppose that $((\pi,W),U)$ acts on $K$. Set $V:=V^{(\pi,W)}$. By definition, $V$ acts on $H=\pi(1_{X_u})K$. 
Since $U:=\{U_{\chi}\}_{\chi \in \widehat{G}}$ commutes with  $\pi(C_0(Y_u))$, it follows that $U_{\chi}$ maps $H$ onto $H$. Clearly, $(U|_{H},V^{(\pi,W)})$ is a weak Weyl pair on $H$
with commuting range projections. We define 
\[
\Phi((\pi,W),U)=(U|_{H},V^{(\pi,W)}).\]

Let $(U,V) \in \mathcal{W}_{c}(P,\widehat{G})$ be given. Suppose that $(U,V)$ acts on $H$. By Thm. \ref{equivalence}, there exists $(\pi,W) \in Rep(C_0(Y_u)\rtimes G)$ such that $V=V^{(\pi,W)}$. Suppose that $(\pi,W)$ acts
on $K$. Then, $H=\pi(1_{X_u})K$. Recall that $W$ is the minimal unitary dilation of $V$. 

Let $\chi \in \widehat{G}$ be given. We claim that there exists a unique unitary operator $\widetilde{U}_\chi$ on $K$ such that 
\begin{enumerate}
\item[(C1)] for $\xi \in H$, $\widetilde{U}_\chi\xi=U_\chi\xi$, and
\item[(C2)] for $x \in G$, $\widetilde{U}_{\chi}W_x=\chi(x)W_{x}\widetilde{U}_\chi$.
\end{enumerate}
Conditions $(C1)$ and $(C2)$ together with  the fact that $\bigcup_{a \in P}W_{a}^{*}H$ is dense in $K$ clearly determine the operator $\widetilde{U}_{\chi}$ uniquely.  

We show the existence below. Define $\widetilde{U}_{\chi}$ on the dense subspace $\bigcup_{a \in P}W_{a}^{*}H$ as follows: for $\xi \in W_{a}^{*}H$, set 
\[
\widetilde{U}_{\chi}\xi=\overline{\chi(a)}W_{a}^{*}U_\chi W_{a}\xi.\]

Let $a,b \in P$ and let $\xi \in W_{a}^{*}H \cap W_{b}^{*}H$ be given. Since $W_xH \subset H$ for $x \in P$, it follows that $W_{a}^{*}H \cap W_{b}^{*}H \subset W_{a+b}^{*}H$. Calculate as follows to observe that 
\begin{align*}
\overline{\chi(a+b)}W_{a+b}^{*}U_{\chi}W_{a+b}\xi&=\overline{\chi(a+b)}W_{a}^{*}W_{b}^{*}U_{\chi}W_{b}W_{a}\xi \\
&=\overline{\chi(a+b)}W_{a}^{*}W_{b}^{*}U_{\chi}V_bW_{a}\xi ~~(~~\textrm{since $W_a\xi \in H$})\\
&=\overline{\chi(a+b)}W_{a}^{*}W_{b}^{*}\chi(b)V_{b}U_{\chi}W_{a}\xi ~~(~~\textrm{since $U_{\chi}V_b=\chi(b)V_{b}U_\chi$})\\
&=\overline{\chi(a)}W_{a}^{*}W_{b}^{*}W_{b}U_{\chi}W_{a}\xi ~~(~~\textrm{since $U_{\chi}W_{a}\xi \in H$})\\
&=\overline{\chi(a)}W_{a}^{*}U_{\chi}W_{a}\xi.
\end{align*}
Similarly, $\overline{\chi(a+b)}W_{a+b}^{*}U_{\chi}W_{a+b}\xi=\overline{\chi(b)}W_{b}^{*}U_{\chi}W_{b}$. This shows that $\widetilde{U}_{\chi}$ is well defined. 
It is clear from the definition that $\widetilde{U}_{\chi}$ is an isometry on $D:=\bigcup_{a \in P}W_{a}^{*}H$ and maps $D$ onto $D$.  Thus, $\widetilde{U}_{\chi}$ extends to a unitary operator to $K$ which 
we again denote by $\widetilde{U}_{\chi}$. 

By definition $\widetilde{U}_{\chi}$ restricted to $H$ coincides with $U_{\chi}$. Using the fact that $(U,V)$ is a weak Weyl pair and the definition of $\widetilde{U}_{\chi}$, it is easy to check, on the dense subspace $D:=\bigcup_{a \in P}W_{a}^{*}H$, that
\[
\widetilde{U}_{\chi}W_{a}=\chi(a)W_{a}\widetilde{U}_\chi\]
for $a \in P$. Since $P$ spans $G$ and $\{W_{x}\}_{x \in G}$ is a group of unitaries, it follows that $\widetilde{U}_{\chi}W_x= \chi(x)W_x\widetilde{U}_{\chi}$ for $x \in G$. Thus, we have established the existence of the unitary operator $\widetilde{U}_{\chi}$ on $K$ for which $(C1)$
and $(C2)$ are satisfied. 

For $x \in G$, let $E_x$ be the projection onto $W_xH$. 
Let $\chi \in \widehat{G}$ be given. By the definition of $\widetilde{U}_{\chi}$, $\widetilde{U}_{\chi}$ maps $W_{a}^{*}H$ onto $W_{a}^{*}H$ for every $a \in P$. Thus, $\widetilde{U}_{\chi}$ commutes with $\{E_{-a}: a \in P\}$. Let $x \in G$ be given. 
Write $x=a-b$ with $a,b \in P$. Then, $E_{x}=W_{a}E_{-b}W_{a}^{*}$. Thanks to the Weyl commutation relation and the fact that $\widetilde{U}_{\chi}$ commutes with $E_{-b}$, it follows that $\widetilde{U}_{\chi}$ commutes with $E_x$ for every $x \in G$.

Let $f \in C_{c}(G)$ be given. Recall that \[
\pi(\widetilde{f})=\int f(x)E_xdx.\]
Since $\widetilde{U}_{\chi}$ commutes with $\{E_x:x \in G\}$,  $\widetilde{U}_{\chi}$ commutes with $\{\pi(\widetilde{f}):f \in C_{c}(G)\}$. Since $\{\widetilde{f}: f \in C_{c}(G)\}$ generates $C_0(Y_u)$, it follows that $\widetilde{U}_{\chi} \in \pi(C_0(Y_u))^{'}$. 

We leave it to the reader to verify that $\widetilde{U}:=\{\widetilde{U}_\chi\}_{\chi \in \widehat{G}}$ is a strongly continuous group of unitaries on $K$. We have now proved that $((\pi,W),\widetilde{U})$ is a representation of the dynamical
system $(C_0(Y_u)\rtimes G,\widehat{G})$. 
Set \[
\Psi(U,V)=((\pi,W),\widetilde{U}).\]

Then, $\Psi$ and $\Phi$ are inverses of each other. We omit this routine verification. \hfill $\Box$

\begin{rmrk}
\label{factorial to factorial}
The maps $\Phi$ and $\Psi$ of Thm. \ref{equivalence1} take factorial representations to factorial representations and take irreducible representations to irreducible representations. The proof that $\Phi$ maps factorial representations to factorial representations proceeds as follows. 

Let $((\pi,W),U)$ be a representation of $(C_0(Y_u)\rtimes G, \widehat{G})$ and let $K$ be the Hilbert space on which it acts. Let $\Phi((\pi,W),U)=(U|_{H},V)$ where $H=\pi(1_{X_u})K$. Denote the von Neumann algebra generated by 
$\{\pi(\phi),W_x,U_{\chi}:\phi \in C_0(Y_u),x \in G, \chi \in \widehat{G}\}$ by $M$ and  the von Neumann algebra on $H$ generated by $\{V_{a},U_{\chi}|_{H}: a \in P, \chi \in \widehat{G}\}$ by $N$. Then, it is routine to prove that 
\[
M^{'} \ni T \to T|_{H} \in N^{'}\]
is an isomorphism between $M^{'}$ and $N^{'}$. Thus, $M$ is a factor if and only if $N$ is a factor. Similarly, $M^{'}=\bbc$ if and only $N^{'}=\bbc$. Thus, $\Phi$ and $\Psi$ map factorial representations to factorial representations and irreducible representations to 
irreducible representations. 
\end{rmrk}

Thm. \ref{main theorem} is an immediate consequence of Thm. \ref{equivalence1} and Takai duality. Recall that Takai duality asserts that $(C_0(Y_u)\rtimes G)\rtimes \widehat{G} \cong C_{0}(Y_u)\otimes \mathcal{K}(L^2(G))$. As a consequence, we have  $Rep((C_0(Y_u)\rtimes G)\rtimes \widehat{G}) \cong Rep(C_0(Y_u))$. The proof of Thm. \ref{main theorem} is essentially transporting the representation theory of $C_0(Y_u)$ to $\mathcal{W}_{c}(P,\widehat{G})$ using Thm. \ref{equivalence1} and by making using of the explicit isomorphism between the $C^{*}$-algebras $(C_0(Y_u) \rtimes G)\rtimes \widehat{G}$ and $C_0(Y_u)\otimes \mathcal{K}(L^2(G))$. For the explicit isomorphism involved in Takai duality, we refer the reader to either \cite{Raeburn_Takai} or \cite{Williams_Dana}. 
We will not write down all the details. For the reader's convenience, we mention a few details concerning  the irreducible weak Weyl pairs with commuting range projections. 

Let us recall the irreducible representations of $(C_{0}(Y_u)\rtimes G)\rtimes \widehat{G}$. Let $K:=L^{2}(G)$. Fix an element $A \in Y_u$. Define a representation $\pi_{A}$ of $C_{0}(Y_u)$ on $K$ by 
\[
\pi_{A}(f)\xi(x)=f(A+x)\xi(x).\]
For $x \in G$, let $W_{x}$ be the unitary on $K$ defined by \[W_{x}\xi(y)=\xi(y-x).\] For $\chi \in \widehat{G}$, let $U_{\chi}$ be the unitary operator on $K$ defined by \[U_{\chi}\xi(y)=\chi(y)\xi(y).\] 
Then, $\{((\pi_A,W),U)\}_{A \in Y_u}$ form a mutually inequivalent exhaustive list of irreducible representations of $(C_0(Y_u)\rtimes G)\rtimes \widehat{G}$. 

For $A \in Y_u$, let $B=-A$ and let  $V^A:=V^{(\pi_A,W)}$. Observe that for $\xi \in K$,
\[
\pi(1_{X_u})\xi(x)=1_{X_u}(A+x)\xi(x)=1_{A}(-x)\xi(x)=1_{B}(x)\xi(x).\]
Thus, $\pi(1_{X_u})K=L^{2}(B)$. By definition, $V^{A}$ is the compression of the left regular representation onto $L^2(B)$.  

Thanks to Thm. \ref{equivalence1}, the assertions in Thm. \ref{main theorem} concerning the irreducible weak Weyl pairs with commuting range projections are now clear. Other assertions can be proved similarly. We leave the details to the reader. 

What about weak Weyl pairs which do not have commuting range projections ? If we drop the assumption that the range projections commute, then we show that, for $P=\bbr_{+}^{2}$,  we can construct weak Weyl pairs that generate 
a factor of both type II and type III. Moreover, we also illustrate that  classifying all the  irreducible weak Weyl pairs is a complicated task. More precisely, we explain a procedure (preserving factoriality and irreducibility) that allows us to build weak Weyl pairs starting from a 
non-degenerate representation of the free product $c_0(\bbn)*c_0(\bbn)$. 

For the rest of this paper, we assume that $P=\bbr_{+}^{2}=[0,\infty)\times [0,\infty)$ and $G=\bbr^2$. We identify $\widehat{G}$ with $\bbr^2$ in the usual way. Let $\{P_{m}\}_{m \geq 1}$ and $\{Q_n\}_{n \geq 1}$ be two sequences of projections on a Hilbert space $K$ such that $P_iP_j=\delta_{ij}P_i$ and $Q_kQ_{\ell}=\delta_{k\ell}Q_k$\footnote{Writing down two such sequences of projections on a Hilbert space $K$ is clearly equivalent to defining a representation of the free product $c_0(\bbn)\ast c_0(\bbn)$.}.

Denote the set of projections on $K$ by $P(K)$. Define a map $F:\bbz^2 \to P(K)$ by 
\begin{equation*}
 F_{(m,n)}:=\begin{cases}
    \sum_{k=1}^{m}P_k & \mbox{if $m \geq 1$ and $ n= 0$}, \cr
      \sum_{k=1}^{n}Q_k & \mbox{if $m=0$ and  $n \geq 1$}, \cr
    1 & \mbox{if~} m \geq 1, n \geq 1, \cr
     0 & \mbox{otherwise}.
         \end{cases}
\end{equation*}
Note that if $(m,n) \in \bbz^2$ and $(p,q) \in \bbn^2$, $F_{(m+p,n+q)} \geq F_{(m,n)}$. 

Let $R:=[0,1] \times [0,1]$ be the unit square and suppose that $\lambda$ is the Lebesgue measure on $R$. Consider $L^{\infty}(R,d\lambda)$ as a $C^{*}$-algebra and let $X$ be the character space of $L^{\infty}(R,d\lambda)$. 
Fix $a,b,c,d \in (0,1)$ such that $a<b$ and $c<d$. Fix a point $z_0 \in X$ such that $1_{[a,b]\times [c,d]}(z_0) \neq 0$. 

Define a map $E:\bbr_{+}^{2} \to P(K)$ as follows. Let $(s,t) \in \bbr_{+}^{2}$ be given. Let $m$ be the integral part of $s$ and let $n$ be the integral part of $t$. 
Set 
\begin{align*}
R_0(s,t):&=[0,m+1-s]\times[0,n+1-t]\\
R_1(s,t):&=[0,m+1-s]\times [n+1-t,1]\\
R_2(s,t):&=[m+1-s,1]\times [0,n+1-t] \\
R_3(s,t):&=[m+1-s,1]\times [n+1-t,1].
\end{align*} 
Define $E_{(s,t)}$ by the following formula. 
\[
E_{(s,t)}:=1_{R_0(s,t)}(z_0)F_{(m,n)}+1_{R_1(s,t)}(z_0)F_{(m,n+1)}+1_{R_2(s,t)}(z_0)F_{(m+1,n)}+1_{R_3(s,t)}(z_0)F_{(m+1,n+1)}.\]
Since $\{1_{R_i}\}_{i=0}^{3}$ is an orthogonal family in $L^{\infty}(R,d\lambda)$ adding up to $1$, exactly one term survives in the above expression. Consequently, $E_{(s,t)}$ is a projection. 

\begin{lmma}
\label{crucial}
With the foregoing notation, we have the following. 
\begin{enumerate}
\item[(1)] The map $E:\bbr_{+}^{2} \to P(K)$ is increasing, i.e  $E_{(s,t)} \leq E_{(s+s_0,t+t_0)}$ for  $(s,t) \in \bbr_{+}^{2}$ and for every $(s_0,t_0) \in \bbr_{+}^{2}$. 
\item[(2)] For $\xi,\eta \in K$, the map 
\[
\bbr_{+}^{2} \ni (s,t) \to \langle E_{(s,t)}\xi|\eta \rangle \in \bbc\]
is Lebesgue measurable. 
\item[(3)] Let $(m,n) \in \bbn^2$ be given. The set $\{(s,t) \in \bbr_{+}^{2}:E_{(s,t)}=F_{(m,n)}\}$ contains a Lebesgue measurable set of positive measure. 
\end{enumerate}
\end{lmma}
\textit{Proof.} The proof of $(1)$ is a case by case verification. Let $(s,t) \in \bbr_{+}^{2}$ be given. Suppose $s_1>s$. Let $m$ be the integral part of $s$, $p$ the integral part of $s_1$ and $n$ the integral part of $t$. 

\textbf{Case 1: $m<p$.}

Let $r:=m+1-s$ and $r_1=p+1-s_1$. 

\textbf{Case (a): $r \leq r_1$.} 

\textbf{Case $(i)$: $1_{R_0(s,t)}(z_0)=1$}. In this case, $E_{(s,t)}=F_{(m,n)}$. Note that $R_{0}(s_1,t)$ contains $R_{0}(s,t)$. Thus, $1_{R_0(s,t)} \leq 1_{R_0(s_1,t)}$ in $L^{\infty}(R,d\lambda)$. Consequently, 
$1_{R_0(s_1,t)}(z_0) =1$. Therefore, $E_{(s_1,t)}=F_{(p,n)}$. Since $F_{(p,n)} \geq F_{(m,n)}$, we have $E_{(s_1,t)}\geq E_{(s,t)}$. 

\textbf{Case $(ii)$: $1_{R_1(s,t)}(z_0)=1$.} We can argue as in Case $(i)$ and deduce $E_{(s_1,t)}\geq E_{(s,t)}$.

\textbf{Case $(iii)$: $1_{R_2(s,t)}(z_0)=1$.} In this case, $E_{(s,t)}=F_{(m+1,n)}$. Note that the union $R_{2}(s_1,t)\cup R_{0}(s_1,t)$ contains $R_{2}(s,t)$. Therefore, either $1_{R_0(s_1,t)}(z_0)=1$ or $1_{R_2(s_1,t)}(z_0)=1$.
This means that $E_{(s_1,t)}$ is either $F_{(p,n)}$ or $F_{(p+1,n)}$. Both $F_{(p,n)}$ and $F_{(p+1,n)}$ are greater than $F_{(m+1,n)}$ as $F$ is increasing and as $p \geq m+1$. Thus, $E_{(s_1,t)} \geq E_{(s,t)}$. 

\textbf{Case $(iv)$: $1_{R_3(s,t)}(z_0)=1$.} In this case, $E_{(s,t)}=F_{(m+1,n+1)}$. Note that the union $R_{3}(s_1,t)\cup R_1(s_1,t)$ contains $R_3(s,t)$. Therefore, either $1_{R_3(s_1,t)}(z_0)=1$ or $1_{R_1(s_1,t)}(z_0)=1$. 
This means that $E_{(s_1,t)}$ is either $F_{(p+1,n+1)}$ or $F_{(p,n+1)}$. In either case, $E_{(s_1,t)} \geq E_{(s,t)}$. 

\textbf{Case (b): $r>r_1$.} The analysis here is similar and we can conclude $E_{(s,t)} \leq E_{(s_1,t)}$. 

\textbf{Case 2: $m=p$.} The analysis here is similar to Case 1 (in this case, Case $(a)$ does not arise) and we can conclude that $E_{(s,t)} \leq E_{(s_1,t)}$. 

Thus, we have proved that $E_{(s,t)} \leq E_{(s+s_0,t)}$ for every $(s,t) \in \bbr_{+}^{2}$ and for every $s_0 \geq 0$. An exactly similar argument shows $E_{(s,t)} \leq E_{(s,t+t_0)}$ for every $(s,t) \in \bbr_{+}^{2}$ and $t_0 \geq 0$. 
Hence, the function $E$ is increasing. This proves $(1)$. 

To prove $(2)$, thanks to the polarisation identity, it suffices to show, that for every $\xi \in K$, the map $\bbr_{+}^{2} \ni (s,t) \to \langle E_{(s,t)}\xi|\xi \rangle \in \bbr$ is Lebesgue measurable. To that effect, let $\xi \in K$ be given and 
define $\phi:\bbr_{+}^{2} \to \bbr$ by 
\[
\phi(s,t):=\langle E_{(s,t)}\xi|\xi \rangle.\]
Then, if we fix one variable, $\phi$ is monotone in the other variable. It is well known (and we leave it to the reader to prove that) that such functions are Lebesgue measurable. This proves $(2)$. 

Let $(m,n) \in \bbn^2$. Let \[A:=\{(s,t) \in [m,m+1)\times [n,n+1): (m+1-s,n+1-t) \in (b,1)\times (d,1)\}.\] Then, $A$ is a Borel set of positive measure. Let $(s,t) \in A$ be given.  Note that 
$R_0(s,t)$ contains $[a,b] \times [c,d]$. Since $1_{[a,b]\times [c,d]}(z_0)=1$, we have $1_{R_0(s,t)}(z_0)=1$. Consequently, for $(s,t) \in A$, $E_{(s,t)}=F_{(m,n)}$. This proves that  
$A \subset \{(s,t) \in \bbr_{+}^{2}: E_{(s,t)}=F_{(m,n)}\}$. The proof of $(3)$ is complete. \hfill $\Box$

Extend $E$ to the whole of $\bbr^2$ by setting $E_{(s,t)}=0$ if $(s,t) \notin \bbr_{+}^{2}$. Then, the extended map $E:\bbr^2 \to P(K)$ is still increasing and Lebesgue measurable. 

Let $L:=L^{2}(\bbr^2,K)$ be the space of square integrable Lebesgue measurable functions taking values in $K$. For $(x,y) \in \bbr^2$, let $U_{(x,y)}$ be the unitary on $L$ defined by 
\[
U_{(x,y)}f(u,v):=e^{i(ux+vy)}f(u,v).\]
For $(s,t) \in \bbr^2$, let $W_{(s,t)}$ be the unitary on $L$ defined by 
\[
W_{(s,t)}f(u,v)=f(u-s,v-t).\]
Define a projection $\widetilde{E}:L \to L$ by 
\[
\widetilde{E}f(u,v)=E_{(u,v)}f(u,v).\]
Set $H:=Ran(\widetilde{E})$. 
Note that $U_{(x,y)}$ commutes with $\widetilde{E}$. Thus, $U_{(x,y)}$ maps $H$ onto $H$. 
We denote the restriction of $U_{(x,y)}$ to $H$ again by $U_{(x,y)}$. 

Using the fact that $E$ is increasing, it is routine to prove that $\widetilde{E}W_{(s,t)}\widetilde{E}=W_{(s,t)}\widetilde{E}$ for every $(s,t) \in \bbr_{+}^{2}$. In other words, $H$ is invariant under $\{W_{(s,t)}:(s,t) \in \bbr_{+}^{2}\}$. For $(s,t) \in \bbr_{+}^{2}$, let $V_{(s,t)}$ be the isometry on $H$ defined by 
\[
V_{(s,t)}=W_{(s,t)}|_{H}.\]

Then, $V:=\{V_{(s,t)}\}_{(s,t) \in \bbr_{+}^{2}}$ is a strongly continuous semigroup of isometries on $H$. Similarly, $U:=\{U_{(x,y)}|_{H}\}_{(x,y) \in \bbr^2}$ is a strongly continuous group of unitaries. Clearly, $(U,V)$ is a weak Weyl pair. 

Let us fix notation. Define 
\begin{align*}
M_0:&=W^{*}\{U_{(x,y)}|_{H}, V_{(s,t)}: (x,y) \in \bbr^2, (s,t) \in \bbr_{+}^{2}\},\\
M_1:&=W^{*}\{U_{(x,y)}, W_{(s,t)},\widetilde{E}: (x,y) \in \bbr^2, (s,t) \in \bbr^2\},\\
N:&=W^{*}\{F_{(m,n)}:(m,n) \in \bbn^2\}=W^{*}\{P_m,Q_n: m \in \bbn, n \in \bbn\}.
\end{align*}
Note that $M_0$ acts on $H$, $M_1$ acts on $L$ and $N$ acts on $K$. 
For a bounded operator $T$ on $K$, let $\widetilde{T}$ be the operator on $L$ defined by 
\[
\widetilde{T}f(y)=Tf(y).\]

\begin{ppsn}
\label{induced}
With the foregoing notation, we have the following. 
\begin{enumerate}
\item[(1)] The map $N^{'} \ni T \to \widetilde{T} \in M_1^{'}$ is an isomorphism. 
\item[(2)] $(W,L)$ is the minimal unitary dilation of $V$. 
\item[(3)] Let $t \in \{I, II, III\}$. The von Neumann algebra  $M_0$ is a factor of type $t$ if and only if $N$ is a factor of type $t$. 
\item[(4)] The weak Weyl pair $(U,V)$ is irreducible if and only if $N^{'}=\bbc$. 
\end{enumerate}
\end{ppsn}
\textit{Proof.} From a routine computation, we see that if $T \in N^{'}$, then $\widetilde{T} \in M_1^{'}$. Let $S \in M_{1}^{'}$ be given. Note that \[W^{*}\{U_{(x,y)}, W_{(s,t)}:(x,y) \in \bbr^2, (s,t) \in \bbr^2\} =B(L^{2}(\bbr))\otimes 1 \subset B(L^{2}(\bbr)\otimes K).\]  Therefore, there exists $T \in B(K)$ such that $S=\widetilde{T}$.  The fact that $S$ commutes with $\widetilde{E}$ translates to the equation
\[
TE_{(s,t)}=E_{(s,t)}T\]
for almost all $(s,t) \in \bbr_{+}^{2}$. By $(3)$ of Lemma \ref{crucial}, $T$ commutes with $F_{(m,n)}$ for every $(m,n) \in \bbn^2$. Thus, $T \in N^{'}$. This completes the proof of $(1)$. 

By definition, $(W,L)$ is a dilation of $V$. Let $Q$ be the projection onto the closure of the subspace $\bigcup_{(s,t) \in \bbr_{+}^{2}}W_{(s,t)}^{*}H$. Note that $Ran(Q)$ is invariant under $U_{(x,y)}$ and $W_{(s,t)}$ for every $(x,y) \in \bbr^2$ and for every $(s,t) \in \bbr^2$. Thus, $Q \in \{U_{(x,y)}, W_{(s,t)}:(x,y),(s,t) \in \bbr^2\}^{'}$. Consequently, $Q=\widetilde{R}$ for some projection $R$ on $K$. 

The condition $Q \geq \widetilde{E}$ translates to the fact that $R \geq E_{(s,t)}$ for almost all $(s,t) \in \bbr_{+}^{2}$. Thanks to $(3)$ of Lemma \ref{crucial}, $R \geq F_{(1,1)}=1$. Thus, $R=1$ and hence $Q=1$. This proves $(2)$. 

As alluded to in Remark \ref{factorial to factorial}, it is not difficult to prove using the fact that $(W,L)$ is the minimal unitary dilation of $V$ that 
\[
M_{1}^{'} \ni T \to T|_{H} \in M_0^{'}\]
is an isomorphism of von Neumann algebras. Now, $(3)$ and $(4)$ follow from $(1)$. \hfill $\Box$

\begin{rmrk}
We conclude this paper with the following remarks. 
\begin{enumerate}
\item[(1)] Thanks to Prop. \ref{induced}, we can construct an irreducible weak Weyl pair starting from an irreducible representation of the free product $c_0(\bbn)\ast c_0(\bbn)$. 
Moreover, inequivalent irreducible representations of $c_0(\bbn)\ast c_0(\bbn)$ lead to inequivalent weak Weyl pairs. Thus, listing out all the irreducible weak Weyl pairs
is at least as hard as describing the dual of $c_0(\bbn)\ast c_0(\bbn)$.  Up to the author's knowledge, a ``good description" of the dual of $c_0(\bbn)\ast c_0(\bbn)$ or even
the dual of some of its natural quotients like $C^{*}(\bbz_n*\bbz_m)$ ($n \geq 2$, $m \geq 3$) is not available in the literature. 

\item[(2)] Observe that for the weak Weyl pair $(U,V)$ constructed in Prop. \ref{induced}, we have $Spec(U)$ is independent of the underlying representation of $c_0(\bbn)\ast c_0(\bbn)$ as long as $F_{(m,n)} \neq 0$ for $(m,n) \in \bbn^2 \backslash \{(0,0)\}$. Thus, Corollary \ref{uniqueness} is not true 
without the commutativity assumption on the range projections. 

\item[(3)] Prop. \ref{induced} allows us to construct weak Weyl pairs that generate a factor of both type II and type III. This is because $c_0(\bbn)\ast c_0(\bbn)$ admit factorial representations of type II and type III
as the $C^{*}$-algebra $c_0(\bbn)\ast c_0(\bbn)$ is not of type I.

\item[(4)]  Let $P$ be a closed convex cone in $\bbr^d$ which
we assume is spanning, i.e. $P-P=\bbr^d$ and pointed, i.e. $P \cap -P=\{0\}$.  Building on the two dimensional case, it is not difficult to construct, in this case, weak Weyl pairs $(U,V)$ that generate a factor of both type II and type III. Also, it is possible to
construct a continuum of irreducible weak Weyl pairs which do not have commuting range projections. 

\end{enumerate}

\end{rmrk}

\bibliography{references}
 \bibliographystyle{amsplain}
 
 \vspace{2.5mm}

\noindent
{\sc S. Sundar}
(\texttt{sundarsobers@gmail.com})\\
{\footnotesize Institute of Mathematical Sciences, \\
A CI of Homi Bhabha National Institute, \\ CIT Campus, 
Taramani, Chennai, 600113, \\ Tamilnadu, India.}\\

\end{document}